\documentclass[twoside,11pt]{article}
\usepackage[preprint]{jmlr2e}
\usepackage{amsmath,mathtools}
\usepackage{times}
\usepackage{booktabs,tabularx,array}
\usepackage{lastpage}
\newcommand{\R}{\mathbb{R}}
\newcommand{\E}{\mathbb{E}}
\newcommand{\M}{\mathcal{M}_{k,d}(R,a,b)}
\newcommand{\op}{\mathrm{op}}
\newcommand{\law}{\mathcal{L}}
\newcommand{\KL}{\operatorname{KL}}
\newcommand{\supp}{\operatorname{supp}}

\numberwithin{equation}{section}

\hypersetup{hidelinks,pdftitle={Sharp Dimension-Dependent Rates for Gaussian Mixtures with an Unknown Common Covariance},pdfauthor={Quan D. Bui}}
\jmlrheading{}{2026}{1--\pageref{LastPage}}{Not submitted}{}{Draft}{Quan D. Bui}
\ShortHeadings{Gaussian Mixtures with Unknown Common Covariance}{Quan D. Bui}
\firstpageno{1}
\begin{document}
\title{Sharp Dimension-Dependent Rates for Gaussian Mixtures\\with an Unknown Common Covariance}
\author{\name Quan D. Bui \email quan.bd@vinuni.edu.vn\\
       \addr Center for AI Research\\
       VinUniversity\\
       Vinhomes Ocean Park, Gia Lam, Hanoi, Vietnam}
\editor{Not assigned (draft)}
\maketitle
\begin{abstract}%
We determine the statistical cost of an arbitrary unknown covariance shared by the components of a finite Gaussian location mixture. For fixed $k\ge2$, we consider mixing measures with at most $k$ atoms in a fixed Euclidean ball and noise covariances whose eigenvalues lie in a fixed positive interval. Uniformly over all sample sizes $n$ and dimensions $d$, the minimax expected Wasserstein loss for the mixing measure is of order $1\wedge(d/n)^{1/(4k)}$, and the minimax expected operator-norm loss for the covariance is of order $1\wedge(d/n)^{1/(2k)}$. The constants are independent of $n,d$, and one jointly feasible measurable estimator attains both rates without logarithmic loss. No separation or minimum mixing weight is required. The upper bound uses established scalar moment comparison and finite-support projection inequalities, together with robust directional moments and a fit over the joint parameter class. For the lower bound, Gaussian quadrature and rank-one covariance compensation produce an exponentially large directional family with identical observation covariance and information of order $t^{4k}$. This identifies a dimension-dependent obstruction that is absent under known or isotropic covariance restrictions. We also give the minimal identifying moment degree and sharp inverse powers. The estimator is computationally unrestricted.
\end{abstract}

\section{Introduction}\label{sec:introduction}

An unknown covariance shared across Gaussian mixture components can hide the direction in which the latent means vary. This changes the global difficulty of recovering the mixing measure, even when the number of components is fixed and every mean lies in a bounded ball. We establish the sharp dependence of that difficulty on the ambient dimension. For $n$ independent observations in $\R^d$ from a mixture with at most $k\ge2$ components, the minimax expected losses are
\[
\begin{aligned}
1\wedge(d/n)^{1/(4k)}&\quad\text{for the mixing measure in }W_1,\\
1\wedge(d/n)^{1/(2k)}&\quad\text{for the covariance in operator norm}.
\end{aligned}
\]
The covariance is arbitrary within a fixed positive spectral interval, and the bounds hold without separation or a lower bound on the mixing weights. All constants are uniform in $n,d$ for fixed model bounds.

The scalar exponents are established. \citet{wu2020} obtained the unknown-common-variance rates $n^{-1/(4k)}$ and $n^{-1/(2k)}$ in one dimension. In higher dimensions, \citet{doss2023} determined the sharp mixing-measure rate under known covariance and gave an upper bound under unknown isotropic covariance. Their discussion preceding Theorem~6.1 explicitly leaves the optimal rate for non-isotropic covariance open. The arbitrary-matrix problem requires identifying a latent direction whose second-moment signature need not be distinguishable from noise.

The obstruction is visible in a two-component family. Put equal latent mass at $tu$ and $-tu$ for an unknown unit direction $u$, and choose the noise covariance to be $vI_d-t^2uu^T$. Every observation law then has covariance $vI_d$. Its first departure from the common Gaussian reference occurs in degree four. More generally, a positive $k$-point Gaussian quadrature law suppresses every discrepancy through degree $2k-1$. There are exponentially many separated directions, but the information to the same Gaussian reference is only of order $t^{4k}$. Comparing $nt^{4k}$ with the logarithm of the packing size gives $t\asymp(d/n)^{1/(4k)}$. The hard family has collinear latent support and only a rank-one change in noise covariance; the phenomenon does not require independently varying every matrix entry.

The matching upper bound estimates the first $2k$ moments in each direction of a fixed sphere net. A common median-of-means partition controls all these estimates, and integrating its tail yields the fractional moment bounds needed for expected parameter loss. We fit one mixing measure and one covariance jointly. Polynomial norming then converts net discrepancies between the fitted and true model moments into uniform directional control. Scalar inverse stability and finite-support transport lifting complete the argument. The construction is measurable and always lies in the original parameter class.

The contribution is the dimension-uniform minimax law for an arbitrary shared covariance, with the covariance rate as a companion conclusion from the same information geometry. The scalar inverse powers are those of \citet[Proposition~3]{wu2020}. The scalar compensated-quadrature construction and its information order already occur in their Proposition~9. The projection mechanism comes from \citet[Lemma~3.1]{doss2023}; using it in the affine direction span of the two finite supports removes ambient dimension from its constant. We give self-contained versions of these ingredients to expose their hypotheses and the uniformity needed here.

Section~\ref{sec:related} positions the result against these and recent neighboring works. Section~\ref{sec:model} states the minimax theorem and an explicit growing-dimension comparison. Section~\ref{sec:inverse} gives inverse stability and its sharpness consequences. Sections~\ref{sec:upper} and~\ref{sec:lower} prove the statistical bounds; the appendices supply all technical arguments.

\section{Related Work and Comparison}\label{sec:related}

\paragraph{Scalar moments and identifiability.}
Moment comparison is central to separation-free finite-mixture estimation. \citet{wu2020} establish polynomial-interpolation inequalities and optimal estimation for Gaussian location mixtures with known or unknown common variance. Their unknown-variance result already has the scalar powers used here, with no minimum-weight condition. Their Gaussian-quadrature lower bound also matches moments through degree $2k-1$ while changing the latent measure and noise variance. We retain those mechanisms and determine how the arbitrary shared matrix changes the dependence on $d$. More broadly, \citet{heinrich2018} characterize local minimax behavior under strong identifiability. Gaussian location-scale families exhibit derivative interactions and weak-identifiability phenomena analyzed by \citet{ho2016}. Those results explain why parameter recovery can be slower than regular estimation, but their hypotheses and local experiments do not give the present dimension-uniform common-matrix risks.

\paragraph{High-dimensional location mixtures.}
The closest rate comparison is \citet{doss2023}. For known identity covariance and bounded mixing support, their Theorem~1.1 gives the sharp $W_1$ rate $(d/n)^{1/4}+n^{-1/(4k-2)}$ in the nonsaturated regime. A fixed known positive-definite covariance with uniformly bounded eigenvalues has the same rate: whitening distorts Euclidean distances and support radii by fixed factors, and inner and outer balls give the corresponding lower and upper comparisons. Their author-manuscript Theorem~6.1 gives the unknown-isotropic upper bound $(d/n)^{1/4}+n^{-1/(4k)}$ and explicitly distinguishes the remaining arbitrary-covariance problem. Table~\ref{tab:comparison} records the comparison with bounded-loss saturation. The theorem and proposition numbers used for the two foundational papers refer to their arXiv author manuscripts, versions 1807.07237v2 and 2002.05818v2, respectively; bibliography entries give the published versions.

\begin{table}[t]
\centering
\small
\begin{tabularx}{\textwidth}{@{}>{\raggedright\arraybackslash}p{0.24\textwidth}p{0.43\textwidth}X@{}}
\toprule
Noise model & Mixing-measure $W_1$ risk or upper bound & Source and scope \\
\midrule
Unknown common variance, $d=1$ & $\asymp n^{-1/(4k)}$ & \citet{wu2020}; sharp \\[3pt]
Known covariance & $\asymp 1\wedge\bigl[(d/n)^{1/4}+n^{-1/(4k-2)}\bigr]$ & \citet{doss2023}; sharp \\[3pt]
Unknown $\sigma^2 I_d$ & $\lesssim 1\wedge\bigl[(d/n)^{1/4}+n^{-1/(4k)}\bigr]$ & \citet{doss2023}; upper bound \\[3pt]
Arbitrary unknown shared $\Sigma$ & $\asymp 1\wedge(d/n)^{1/(4k)}$ & Theorem~\ref{thm:minimax}; sharp \\
\bottomrule
\end{tabularx}
\caption{Theoretical comparison for a fixed upper bound $k\ge2$ on the number of atoms, bounded Euclidean mixing support, and fixed positive noise spectral bounds. Constants can depend on those fixed bounds. The known-covariance row follows by whitening. The isotropic row records the published upper bound only. All losses are expected Euclidean $W_1$ losses; saturation follows from bounded mixing support.}
\label{tab:comparison}
\end{table}

The isotropic row is used as an upper bound. A scalar unknown-variance lower bound cannot simply be embedded unchanged into every isotropic dimension: when $d\ge k$, the centered covariance of an at-most-$k$-atomic latent law has rank at most $k-1$, so the smallest eigenvalue of the population observation covariance equals the scalar noise variance. In contrast, the compensated family used below varies the noise covariance with $u$ and is outside the isotropic submodel for $d>1$.

\paragraph{Shared scale and Gaussian cancellation.}
\citet{bariletto2026} study infinite homoscedastic location-scale mixtures, including latent transport and operator-norm scale discrepancies. Their Gaussian inverse bounds and posterior consequences allow infinite mixing measures and have dimension-dependent constants and logarithmic or iterated-logarithmic inverse behavior. A fixed finite support bound changes the inverse modulus and permits the polynomial dimension law studied here. \citet{zhang2026} establish the exact scalar cancellation order in the Ho--Nguyen polynomial system using a common variance shift and positive Gauss--Hermite quadrature, with local statistical consequences. The cancellation order and quadrature witness are therefore established mathematical background; the result here concerns uniform global risks as $d$ varies under one arbitrary shared covariance.

\paragraph{Clustering, algorithms, and density loss.}
Unknown shared covariance has also been studied for clustering. \citet{davis2025} analyze misclassification under an unknown, potentially ill-conditioned common covariance. \citet{bing2025} study multicomponent EM with an unknown full covariance under separation and initialization conditions. These works study clustering or separated component recovery. Our global latent $W_1$ class includes collisions and vanishing weights. For known unit covariance, \citet{he2026} obtain sharp proper density estimation in polynomial arithmetic time for fixed $k$, measured in Hellinger distance. Density loss and latent-parameter loss capture different difficulties. The theorem below specifies the statistical limit for the latter with covariance uncertainty and an unrestricted estimator.

\section{Model and Sharp Minimax Rates}\label{sec:model}

Fix an integer $k\ge2$, a radius $R>0$, and $0<a<b<\infty$. Let $\overline B_R\subset\R^d$ be the closed Euclidean ball centered at zero, and define
\[
\M=\left\{(P,\Sigma): |\supp P|\le k,\ \supp P\subseteq\overline B_R,
\quad aI_d\preceq\Sigma\preceq bI_d\right\},
\]
where $P$ is a probability measure and $\Sigma$ is symmetric. The experiment observes $X_1,\ldots,X_n$ independently with common law
\[
F_{P,\Sigma}=P*N(0,\Sigma).
\]
Equivalently, $X=Z+\varepsilon$ with independent $Z\sim P$ and $\varepsilon\sim N(0,\Sigma)$. Write $E$ or $\E$ for expectation, with a subscript when the parameter needs to be specified.

For probability measures with finite first moments, use Euclidean transport cost:
\[
W_1(P,Q)=\inf_{\pi\in\Pi(P,Q)}\int\|x-y\|\,d\pi(x,y).
\]
For a matrix $A$, $\|A\|_{\op}=\sup_{\|u\|=1}\|Au\|$; for symmetric $A$ this equals $\sup_{\|u\|=1}|u^TAu|$. Define the two minimax risks by
\begin{align*}
\mathcal R_P(n,d)&=\inf_{\widehat P}\ \sup_{(P,\Sigma)\in\M}
  \E_{P,\Sigma}W_1(\widehat P,P),\\
\mathcal R_\Sigma(n,d)&=\inf_{\widehat\Sigma}\ \sup_{(P,\Sigma)\in\M}
  \E_{P,\Sigma}\|\widehat\Sigma-\Sigma\|_{\op}.
\end{align*}
The infima range over Borel estimators taking values in finite-first-moment probability measures and real $d\times d$ matrices, respectively. The lower bounds also apply to probability-valued outputs with infinite loss. The feasible estimator used for the upper bounds is a special member of these unrestricted estimator classes. Labels of mixture components play no role.

\begin{theorem}[Sharp minimax rates]\label{thm:minimax}
For every fixed $k\ge2,R>0,0<a<b<\infty$, there exist constants $0<c<C<\infty$, depending only on $k,R,a,b$, such that for every pair of positive integers $n,d$,
\begin{equation}\label{eq:source-1-1}
c\min\{1,(d/n)^{1/(4k)}\}\le\mathcal R_P(n,d)
\le C\min\{1,(d/n)^{1/(4k)}\},
\end{equation}
\begin{equation}\label{eq:source-1-2}
c\min\{1,(d/n)^{1/(2k)}\}\le\mathcal R_\Sigma(n,d)
\le C\min\{1,(d/n)^{1/(2k)}\}.
\end{equation}
One Borel measurable estimator $(\widehat P,\widehat\Sigma)\in\M$ attains both upper bounds. Neither the truth nor the fit requires separated atoms or a positive weight floor.
\end{theorem}

The proof of the upper bounds is in Section~\ref{sec:upper} and Appendix~\ref{app:estimator}; the matching lower bounds are in Section~\ref{sec:lower} and Appendix~\ref{app:testing}. Both losses are unsquared. In particular, the covariance in the theorem is the Gaussian noise covariance, not the total covariance of the observations.

Uniform consistency along a sequence $d=d_n$ holds precisely when $d_n/n\to0$. If the ratio stays bounded away from zero along a subsequence, the lower bounds stay positive. The constants may grow with $k$, $R$ or the spectral bounds. Thus fixing $k$ before varying $n,d$ is part of the theorem. The fit is an information-theoretic construction with unrestricted computation.

\paragraph{A concrete dimension effect.}
For $k=2$ and $d\asymp\sqrt n$, Theorem~\ref{thm:minimax} gives mixing-measure risk $\asymp n^{-1/16}$ and covariance operator risk $\asymp n^{-1/8}$. Under known covariance, the mixing-measure rate in Table~\ref{tab:comparison} becomes $n^{-1/8}+n^{-1/6}\asymp n^{-1/8}$. The published unknown-isotropic upper bound also becomes $O(n^{-1/8})$. Therefore allowing an arbitrary shared covariance causes a strict deterioration in global latent recovery relative to the isotropic restriction. This comparison uses only its upper bound. Table~\ref{tab:dimension} makes the substitutions explicit.

\begin{table}[t]
\centering
\begin{tabularx}{\textwidth}{@{}Xcc@{}}
\toprule
Noise model & Substituted $W_1$ expression & Order \\
\midrule
Known covariance & $n^{-1/8}+n^{-1/6}$ & $n^{-1/8}$ \\[3pt]
Unknown isotropic (upper bound) & $n^{-1/8}+n^{-1/8}$ & $O(n^{-1/8})$ \\[3pt]
Arbitrary unknown shared covariance & $(n^{1/2}/n)^{1/8}$ & $n^{-1/16}$ \\
\bottomrule
\end{tabularx}
\caption{An algebraic dimension-effect illustration: $k=2$ and $d\asymp\sqrt n$ as $n\to\infty$, with fixed $R,a,b$. The first two rows substitute into Table~\ref{tab:comparison}; the last row is the sharp rate of Theorem~\ref{thm:minimax}. The isotropic row is an upper bound.}
\label{tab:dimension}
\end{table}

\section{Directional Moments and Inverse Stability}\label{sec:inverse}

For two observation laws in the model, define the directional raw-moment discrepancy
\[
D_{2k}(F_{P,\Sigma},F_{Q,\Lambda})
=\max_{1\le r\le2k}\ \sup_{\|u\|=1}
\left|\E_{P,\Sigma}(u^TX)^r-\E_{Q,\Lambda}(u^TX)^r\right|.
\]
The bounded support and noise spectral bounds make all these moments finite. Each directional moment is a homogeneous polynomial in $u$. Consequently, equality of directional moments through degree $q$ is equivalent to equality of all multivariate raw moments of total degree at most $q$.

\begin{theorem}[Uniform inverse stability]\label{thm:inverse}
Under the fixed model bounds of Theorem~\ref{thm:minimax}, there is a constant $C=C(k,R,a,b)$ such that, for every $d\ge1$ and every $(P,\Sigma),(Q,\Lambda)\in\M$, writing $D=D_{2k}(F_{P,\Sigma},F_{Q,\Lambda})$,
\begin{equation}\label{eq:source-1-3}
W_1(P,Q)\le C D^{1/(2k)},\qquad
\|\Sigma-\Lambda\|_{\op}\le C D^{1/k}.
\end{equation}
\end{theorem}

\begin{corollary}[Minimal moment degree and sharp inverse powers]\label{cor:sharp}
For each $d\ge1$ in the same model class, the directional raw moments through degree $2k$ identify $(P,\Sigma)$. Moments through degree $2k-1$ do not identify either parameter uniformly over the class. Moreover, no exponent larger than $1/(2k)$ in the first inequality of \eqref{eq:source-1-3}, or larger than $1/k$ in the second, can hold uniformly, even with $d$ fixed.
\end{corollary}

The identification assertion follows immediately from Theorem~\ref{thm:inverse}. The converse and sharpness assertions are proved using the admissible quadrature family in Section~\ref{sec:sharpness}. This is a statement about the largest available moment degree, not a lower bound on the number of arbitrary identifying statistics.

The scalar ingredient has the same powers as \citet[Proposition~3]{wu2020}. We state it with variances as its scale parameters.

\begin{lemma}[Scalar moment comparison]\label{lem:scalar}
Let $P,Q$ be probability measures with at most $k$ atoms in $[-R,R]$ and let $s,t\in[0,b]$ be variances. Set
\[
\delta=\max_{1\le r\le2k}
\left|m_r(P*N(0,s))-m_r(Q*N(0,t))\right|,
\qquad m_r(A)=\int x^r\,dA(x).
\]
There is a constant depending only on $k,R,b$ such that
\[
|s-t|\le C\delta^{1/k},\qquad W_1(P,Q)\le C\delta^{1/(2k)}.
\]
The conclusions include $\delta=0$, fewer than $k$ distinct atoms, and vanishing weights.
\end{lemma}

Appendix~\ref{app:scalar} proves the lemma directly by deconvolving the smaller variance and comparing the resulting law with a finite measure. Gap-independent Hermite interpolation controls transport; the squared annihilating polynomial and Gaussian Hermite orthogonality control the variance difference. This also verifies the normalization and all degeneracies needed when applying the cited scalar result.

\begin{lemma}[Finite-support transport projection]\label{lem:projection}
Let $P$ and $Q$ each be a probability measure on at most $k$ atoms in $\R^d$. There is a unit direction $u$ such that
\[
W_1(P,Q)\le c_k^{-1}W_1(u_\#P,u_\#Q),\qquad
c_k=\left[16\binom{2k}{2}\sqrt{2k-1}\right]^{-1}.
\]
Here $u_\#P$ is the law of $u^TZ$ for $Z\sim P$.
\end{lemma}

\begin{proof}
Let \(S=\operatorname{supp}P\cup\operatorname{supp}Q\), and \(m=|S|\le2k\). If \(m=1\), then \(P=Q\) and any unit direction suffices. Otherwise let

\[
H=\operatorname{span}\{x-y:x,y\in S\},\quad
s=\dim H\le m-1,\quad N={m\choose2}.
\]

Here \(H\) is the direction space of the affine hull, even if the affine hull does not contain the origin. Let \(g\) be standard Gaussian on \(H\). By Markov,
\(\Pr(\|g\|>2\sqrt s)\le1/4\). For each normalized nonzero support difference \(v\),

\[
\Pr(|g^Tv|<1/(8N))\le\sqrt{2/\pi}/(8N)\le1/(8N).
\]

The union bound over all \(N\) differences shows that some \(g\) satisfies both \(\|g\|\le2\sqrt s\) and all \(|g^Tv|\ge1/(8N)\). Normalize it to a unit \(u\in H\). Then

\begin{equation}\label{eq:source-3-1}
|u^T(x-y)|\ge c_k\|x-y\|,
\qquad c_k=\left[16{2k\choose2}\sqrt{2k-1}\right]^{-1}
\end{equation}

for every \(x,y\in S\).

Projection by \(u\) is injective on \(S\). An optimal coupling of the two projected finite measures is a matrix indexed by projected support points. Injectivity identifies every row and column with its original atom, so the same matrix lifts to a coupling of \(P,Q\). Inequality \eqref{eq:source-3-1} bounds each Euclidean cost by \(c_k^{-1}\) times its projected cost. Therefore

\begin{equation}\label{eq:source-3-2}
W_1(P,Q)\le c_k^{-1}W_1(u_\#P,u_\#Q).
\end{equation}

This is the finite-support mechanism in \citet{doss2023}, Lemma 3.1. Its stated ambient-dimension factor can be replaced here by a factor depending only on \(k\), by working in \(H\), whose dimension is at most \(2k-1\).

\end{proof}

\subsection{Proof of Theorem~\ref{thm:inverse}}
\begin{proof}
For every ambient unit \(u\), the two projected mixing laws have at most \(k\) points in \([-R,R]\) and scalar noise variances \(u^T\Sigma u,u^T\Lambda u\in[a,b]\). Their moment discrepancy is bounded by \(D\). Apply Lemma~\ref{lem:scalar}:

\begin{equation}\label{eq:source-3-3}
W_1(u_\#P,u_\#Q)\le C D^{1/(2k)},\qquad
|u^T(\Sigma-\Lambda)u|\le C D^{1/k}.
\end{equation}

For the first assertion choose the direction in \eqref{eq:source-3-1} and apply \eqref{eq:source-3-2}. Choosing this direction after seeing \(P,Q\) is valid because \(D\) already controls every direction. For the second assertion take the supremum over \textbf{all} ambient unit directions. The Rayleigh-quotient formula for a symmetric matrix gives

\[
\|\Sigma-\Lambda\|_{\rm op}
=\sup_{\|u\|=1}|u^T(\Sigma-\Lambda)u|\le C D^{1/k}.
\]

This includes directions outside the affine support span. These arguments prove \eqref{eq:source-1-3}. If the observation laws agree, all their moments agree, so \eqref{eq:source-1-3} proves identifiability of the mixing measure and covariance, with labels irrelevant.
\end{proof}

\section{A Joint Estimator without Logarithmic Loss}\label{sec:upper}

The inverse theorem reduces estimation to uniform directional moment control. The index set has a sphere of directions rather than all entries of a moment tensor. Its exponential net size contributes order $d$ to the block count, while the variance of each directional statistic remains bounded independently of $d$.

Put $q=2k$, $c_r=r^r/r!$, $L_q=\max_{1\le r\le q}rc_r$, and $\zeta=(2L_q)^{-1}$. Fix a deterministic $\zeta$-net $\mathcal N_d$ of the unit sphere with $|\mathcal N_d|\le(1+2/\zeta)^d$. Polarization gives, for every real homogeneous polynomial $h$ of degree $1\le r\le q$,
\begin{equation}\label{eq:norming-main}
\sup_{\|u\|=1}|h(u)|\le2\max_{v\in\mathcal N_d}|h(v)|.
\end{equation}
The proof, including a constant independent of ambient dimension, is in Appendix~\ref{app:estimator}.

Choose $B$ to be the smallest odd integer satisfying
\begin{equation}\label{eq:blocks-main}
B\ge3,\qquad B\ge\log_2q+d\log_2(1+2/\zeta).
\end{equation}
Then $B\asymp_k d$. If $n\ge2B$, let $m=\lfloor n/B\rfloor$ and divide the first $Bm$ observations into $B$ deterministic disjoint blocks of size $m$. For every $(r,v)\in\{1,\ldots,q\}\times\mathcal N_d$, take the median $M_{r,v}$ of the block averages of $(v^TX)^r$. All coordinates use the same blocks. Median-of-means estimation is a standard finite-variance device; see, for example, \citet{devroye2016}. We prove the simultaneous bound needed here directly.

Let $V\ge1$ be a common bound on $\E|v^TX|^{2r}$, supplied by \eqref{eq:source-4-1}, and write
\[
\mathcal E=\max_{r,v}|M_{r,v}-\E_{P,\Sigma}(v^TX)^r|,
\qquad t_0=8\sqrt{V/m}.
\]
Chebyshev's inequality within a block and independence across blocks imply
\begin{equation}\label{eq:tail-main}
\Pr(\mathcal E>t_0s)\le2^{-B}s^{-B}\quad(s\ge1),
\qquad
\E\mathcal E^p\le C_{k,R,b,p}(d/n)^{p/2}\quad(0<p\le1).
\end{equation}
The second inequality integrates the first for this same fixed partition. Independence between coordinates is unnecessary.

We fit the targets over one compact parameter space,
\[
\Theta_d=\Delta_k\times(\overline B_R)^k
\times\{S=S^T:aI_d\preceq S\preceq bI_d\},
\]
allowing repeated locations and zero weights. For $\theta=(w,\mu,S)$, its model moments $m_{r,v}(\theta)$ are the continuous polynomials in \eqref{eq:source-5-2}. Fix a countable dense sequence $(\theta_j)_{j\ge1}$ and minimize the maximum target error over that sequence up to tolerance $\tau_n=\sqrt{B/n}$: choose the first $j$ within $\tau_n$ of the infimum. Formula~\eqref{eq:source-5-3} specifies the estimator exactly. The dense infimum equals the full feasible minimum. The first admissible index is Borel measurable because each objective and the countable infimum are Borel. Appendix~\ref{app:estimator} verifies density even at boundary parameters and gives the selection events explicitly.

\begin{proposition}[Moment risk of the joint fit]\label{prop:upper}
For the block count \eqref{eq:blocks-main}, if $n\ge2B$, the preceding estimator is Borel measurable, takes values in $\M$, and satisfies, for every $0<p\le1$,
\[
\sup_{(P,\Sigma)\in\M}
\E_{P,\Sigma}\!\left[
D_{2k}(F_{\widehat P,\widehat\Sigma},F_{P,\Sigma})^p\right]
\le C_{k,R,b,p}(d/n)^{p/2}.
\]
The estimator and block partition are the same for every $p$.
\end{proposition}

\begin{proof}
The true parameter is feasible, so the fitted net error is at most $\mathcal E+\tau_n$. Hence the net moment difference between the fitted and true models is at most $2\mathcal E+\tau_n$. For each $r$, this difference is a genuine homogeneous polynomial in the direction. Apply \eqref{eq:norming-main} to obtain
\[
D_{2k}(F_{\widehat P,\widehat\Sigma},F_{P,\Sigma})\le4\mathcal E+2\tau_n.
\]
Subadditivity of the $p$th power for $0<p\le1$, \eqref{eq:tail-main}, and $B\le C_kd$ prove the claim. The median targets themselves need not be a polynomial; norming is used only after the joint fit. The full construction and quantitative tail integration are in Appendix~\ref{app:estimator}.
\end{proof}

\paragraph{Completion of the upper bounds.}
Apply Theorem~\ref{thm:inverse} and Proposition~\ref{prop:upper} with $p=1/(2k)$ and $p=1/k$. This gives the two upper rates for $n\ge2B$. For $n<2B$, use the constant feasible estimator $(\delta_0,aI_d)$, whose losses are at most $R$ and $b-a$. Since $n<2C_kd$, each saturated target in Theorem~\ref{thm:minimax} is bounded below by a positive fixed constant. These bounded losses therefore prove the upper bounds for every $n,d$. Appendix~\ref{app:estimator} includes the explicit inequalities for both branches.

\section{Covariance Compensation and Matching Lower Bounds}\label{sec:lower}

The lower bound lifts the scalar quadrature alternative of \citet[Proposition~9]{wu2020} to a large directional family. Adjusting a shared matrix rather than a scalar variance lets every direction have the same total observation covariance. The resulting family connects the scalar information order to ambient dimension.

\subsection{Two Symmetric Components}

Fix $v=(a+b)/2$, a unit vector $u$, and $t>0$ with $t\le R$, $t^2\le(v-a)/2$ and $t\le\sqrt v/2$. Define
\begin{equation}\label{eq:source-6-1}
P_{u,t}=\tfrac12\delta_{tu}+\tfrac12\delta_{-tu},\qquad
\Sigma_{u,t}=vI_d-t^2uu^T,\qquad Q_0=N(0,vI_d).
\end{equation}
The noise eigenvalues are $v-t^2$ along $u$ and $v$ on its orthogonal complement, so the pair is admissible. The observation law $F_{u,t}=F_{P_{u,t},\Sigma_{u,t}}$ has a positive normalized density. Writing $y=u^Tx$, direct division of its two Gaussian densities by the reference density gives
\begin{equation}\label{eq:source-6-2}
\begin{aligned}
\frac{dF_{u,t}}{dQ_0}(x)
={}&\sqrt{\frac v{v-t^2}}
\exp\!\left\{-\frac{t^2}{2(v-t^2)}
-\frac{t^2y^2}{2v(v-t^2)}\right\}
\cosh\!\left(\frac{ty}{v-t^2}\right).
\end{aligned}
\end{equation}
Every law in this family has mean zero and covariance $vI_d$. Its first three directional raw moments agree with the reference, and for every $z\in\R^d$,
\begin{equation}\label{eq:source-6-3}
\E(z^TX)^{4}-\E_{Q_0}(z^TX)^{4}=-2t^4(z^Tu)^4.
\end{equation}
Indeed, the projected latent variable has second moment $t^2(z^Tu)^2$ and fourth cumulant $-2t^4(z^Tu)^4$, while the noise contributes no fourth cumulant. The general construction below bounds the chi-square divergence by $Ct^8$, with an integrable remainder.

\subsection{A General Compensated Family}

\begin{proposition}[Quadrature alternatives and information]\label{prop:alternatives}
For the monic probabilists' Hermite polynomials $H_j$, there is a symmetric random variable $A$ supported on the $k$ distinct roots of $H_k$. Its weights are fixed and strictly positive, and for $G\sim N(0,1)$,
\[
\E A^j=\E G^j\quad(0\le j\le2k-1),\qquad
\E H_{2k}(A)=-k!.
\]
Let $B_k=\max|\supp A|$ and set
\[
t_{\max}=\min\{1,R/B_k,\sqrt{(v-a)/2},\sqrt v/2\}.
\]
For every $d\ge1$, unit $u\in\R^d$ and $0<t\le t_{\max}$, the pair
\[
P_{u,t}=\law(tAu),\qquad\Sigma_{u,t}=vI_d-t^2uu^T
\]
is in $\M$. All multivariate raw moments of its observation law through total degree $2k-1$ agree with $Q_0$, and
\begin{equation}\label{eq:moment-main}
\E_{u,t}(z^TX)^{2k}-\E_{Q_0}(z^TX)^{2k}
=-k!t^{2k}(z^Tu)^{2k}\qquad(z\in\R^d).
\end{equation}
There are constants $c,C>0$ depending only on $k,a,b$ such that, uniformly in the stated $d,u,t$,
\[
c\,t^{4k}\le\chi^2(F_{u,t},Q_0)\le C t^{4k},
\qquad \KL(F_{u,t}\|Q_0)\le C t^{4k}.
\]
\end{proposition}

\begin{proof}
Appendix~\ref{app:quadrature} constructs the positive quadrature weights and proves every assertion. The key information calculation can be seen by scaling the projection along $u$ to $Y_s=sA+\sqrt{1-s^2}G$, where $s=t/\sqrt v\le1/2$. Its likelihood ratio $L_s$ against the standard normal density $\phi$ obeys
\[
0<L_s(x)\le \frac2{\sqrt3}\exp(B_k^2/2)\quad(x\in\R,\ 0\le s\le1/2).
\]
Thus $L_s$ is uniformly in $L^2(\phi)$, with an explicit Gaussian tail bound. Hermite completeness and Parseval give the exact nonnegative series
\[
\chi^2(Y_s,G)=\sum_{m\ge2k}\frac{s^{2m}(\E H_m(A))^2}{m!}.
\]
The $m=2k$ term gives the lower bound. Comparing each summand to its value at $s=1/2$ gives the upper bound of order $s^{4k}$. This comparison bounds the full nonnegative remainder. The common independent orthogonal Gaussian factor does not affect the divergence, and $\KL\le\chi^2$.
\end{proof}

\subsection{Directional Testing}

Let $\rho=u^Tu'$. Every atom of $P_{u',t}$ lies on the line spanned by $u'$. The distance of each source atom to that line gives
\begin{equation}\label{eq:distances-main}
\begin{aligned}
W_1(P_{u,t},P_{u',t})&\ge t\E|A|\sqrt{1-\rho^2}
\ge (t/B_k)\sqrt{1-\rho^2},\\
\|\Sigma_{u,t}-\Sigma_{u',t}\|_{\op}&=t^2\sqrt{1-\rho^2}.
\end{aligned}
\end{equation}
The covariance identity follows from the two nonzero eigenvalues of $uu^T-u'u'^T$. Appendix~\ref{app:testing} supplies the calculations, including zero quadrature nodes and collinear directions.

For all $d\ge d_0$, a projective packing gives $M\ge\exp(c_0d)$ directions satisfying $|u_i^Tu_j|\le1/2$ for $i\ne j$. Packing modulo sign is necessary because $A$ is symmetric. Choose a fixed sufficiently small $\eta_0>0$ and
\[
t=\eta_0\min\{1,(d/n)^{1/(4k)}\}.
\]
By Proposition~\ref{prop:alternatives}, the average product KL to $Q_0^{\otimes n}$ is at most $nCt^{4k}\le(1/8)\log M$. Reference-law Fano testing has error bounded below by a positive constant. Nearest-member decoding and \eqref{eq:distances-main} give expected losses at least constant multiples of $t$ and $t^2$, separately for arbitrary estimators under the two losses. When $n\le d$, the same argument uses a fixed positive $t$ and yields bounded-loss saturation.

For $1\le d<d_0$, compare $(\delta_0,vI_d)$ with $(P_{e_1,t},vI_d-t^2e_1e_1^T)$, at $t=\eta_1n^{-1/(4k)}$ with fixed small $\eta_1>0$. Their product KL is at most $1/8$, and their parameter distances are $t\E|A|$ and $t^2$. Pinsker's inequality and two-point testing yield $cn^{-1/(4k)}$ and $cn^{-1/(2k)}$. Because $d$ is bounded, these dominate the saturated targets up to fixed constants. Appendix~\ref{app:testing} gives the packing, testing inequalities and constants for all regimes. This completes the lower bounds of Theorem~\ref{thm:minimax}.

For $d\ge d_0$, the lower bound already holds on the subclass with the fixed quadrature weights, exactly $k$ collinear atoms and rank-one covariance compensation. The weights are strictly positive independently of $n,d,t$. The obstruction in that family is covariance confounding of a direction; its atoms can approach one another as $t$ decreases. The bounded-dimension proof uses the additional one-atom reference and is a statement for the at-most-$k$ class.

\subsection{Proof of Corollary~\ref{cor:sharp}}\label{sec:sharpness}
\begin{proof}
Fix any $d\ge1$ and unit $u$. Both $(P_{u,t},\Sigma_{u,t})$ from Proposition~\ref{prop:alternatives} and $(\delta_0,vI_d)$ are admissible for all $0<t\le t_{\max}$. Their moments agree through degree $2k-1$, but $P_{u,t}\ne\delta_0$ and $\Sigma_{u,t}\ne vI_d$. Thus that degree does not identify either parameter. At degree $2k$, \eqref{eq:moment-main} gives exactly
\[
D_{2k}(F_{u,t},Q_0)=k!t^{2k},\qquad
W_1(P_{u,t},\delta_0)=t\E|A|,\qquad
\|\Sigma_{u,t}-vI_d\|_{\op}=t^2.
\]
Here $\E|A|\ge B_k^{-1}>0$. If a uniform transport inverse held with exponent $\beta>1/(2k)$, these identities would require $t\E|A|\le C(k!)^\beta t^{2k\beta}$ for arbitrarily small $t$, which is impossible. The covariance argument is the same with $t^2$ and $\beta>1/k$. Sufficiency at degree $2k$ follows from Theorem~\ref{thm:inverse} at $D=0$.
\end{proof}

\section{Discussion}\label{sec:discussion}

An arbitrary shared Gaussian covariance changes the dimension dependence of global latent recovery. The hard family has a common total covariance, so the second moment does not reveal its latent direction. Positive quadrature postpones the first discrepancy to degree $2k$, and a large directional packing makes the corresponding scalar information order interact with ambient dimension. The matching joint moment estimator shows that this obstruction determines both the mixing-measure and covariance risks.

\paragraph{Scope and limitations.}
The results hold for a fixed known upper bound $k\ge2$, compact Euclidean mixing support, and a fixed positive spectral interval with $a<b$. Constants depend on these quantities. The model requires one shared covariance, and the losses are expected $W_1$ and covariance operator norm. The fit can search exponentially many directions and a countable dense parameter set; the theorem provides a measurable existence construction with no computational complexity guarantee. Adaptation to growing or unknown $k$, component-specific covariances, and other covariance losses requires additional analysis.

Structured covariance restrictions offer a route to understanding the transition between isotropic and arbitrary-matrix noise. Efficient procedures that attain the unrestricted statistical limits pose a separate question. The rates proved here give a precise benchmark for those investigations: the cost of hiding a latent direction is already present with fixed positive weights, collinear means and rank-one covariance compensation.

\acks{This draft leaves the author list, affiliations, funding declaration, competing-interest declaration, and any required assistance disclosure for completion by the authors before submission.}
\appendix
\section{Scalar Moment Comparison}\label{app:scalar}

Write \(m_r(A)=\int x^r\,dA(x)\).

\subsection{A comparison lemma with only one finite law}

Suppose \(Q\) is a probability measure with distinct support points \(q_1<\cdots<q_\ell\) in \([-R,R]\), and \(A\) is any probability measure with finite absolute moment of order \(2\ell\). If

\[
\max_{1\le r\le2\ell}|m_r(A)-m_r(Q)|\le\eta,\qquad 0<\eta\le1,
\]

then

\begin{equation}\label{eq:source-2-1}
W_1(A,Q)\le C_{\ell,R}\eta^{1/(2\ell)}.
\end{equation}

If all these moments agree, the Wasserstein distance is zero.

\begin{proof} Take any 1-Lipschitz function \(f\) with \(f(0)=0\). Fix a smooth nonnegative function \(\rho\) supported on \([-1,1]\) with integral one, and put \(f_\epsilon=f*\rho_\epsilon\), \(0<\epsilon\le1\). Cancellation of the derivative of a compactly supported mollifier and the Lipschitz inequality imply

\begin{equation}\label{eq:source-2-2}
\|f-f_\epsilon\|_\infty\le\epsilon,
\qquad \|f_\epsilon^{(j)}\|_\infty\le C_j\epsilon^{1-j}
\quad(j\ge1).
\end{equation}

For example, subtract \(f(x)\) in the integral for the \(j\)-th derivative and use
\(\int |y|\,|\rho_\epsilon^{(j)}(y)|\,dy=C_j\epsilon^{1-j}\).

Let \(p_\epsilon\), of degree at most \(2\ell-1\), interpolate both \(f_\epsilon\) and its first derivative at every \(q_i\). Existence and uniqueness follow because the homogeneous interpolation kernel would have \(\ell\) double roots and degree less than \(2\ell\). In Newton form its nodes are the repeated list
\((q_1,q_1,\ldots,q_\ell,q_\ell)\). A divided difference of order \(j\ge1\) is bounded by \(\|f_\epsilon^{(j)}\|_\infty/j!\), by repeated Rolle's theorem, including coincident nodes by continuity. The order-zero coefficient has absolute value at most \(R+1\). The monomial coefficient sum of a Newton product of length \(j\) is at most \((1+R)^j\). It follows that

\begin{equation}\label{eq:source-2-3}
\|p_\epsilon\|_{\mathrm{coef},1}
\le C_{\ell,R}\epsilon^{2-2\ell}.
\end{equation}

Here \(\|\cdot\|_{\mathrm{coef},1}\) is the sum of absolute monomial coefficients. This bound is independent of all gaps between support points.

Put \(q(x)=\prod_{i=1}^\ell(x-q_i)^2\). The Hermite interpolation remainder, by Rolle's theorem on the interval containing \(x\) and all the nodes, gives for \textbf{every} real \(x\)

\begin{equation}\label{eq:source-2-4}
|f_\epsilon(x)-p_\epsilon(x)|
\le \frac{\|f_\epsilon^{(2\ell)}\|_\infty}{(2\ell)!}q(x)
\le C_\ell\epsilon^{1-2\ell}q(x).
\end{equation}

It is zero at a node. As \(E_Qq=0\), \(q\ge0\), and
\(\|q\|_{\mathrm{coef},1}\le(1+R)^{2\ell}\), the moment assumption yields

\begin{equation}\label{eq:source-2-5}
0\le E_Aq\le(1+R)^{2\ell}\eta.
\end{equation}

The absolute-moment hypothesis ensures finite expectations. Combining \eqref{eq:source-2-2}–\eqref{eq:source-2-5}, with equality of \(p_\epsilon,f_\epsilon\) on the support of \(Q\) gives

\[
\begin{aligned}
|E_Af-E_Qf|
&\le2\epsilon+|E_Ap_\epsilon-E_Qp_\epsilon|
 +E_A|f_\epsilon-p_\epsilon|\\
&\le C_{\ell,R}
 (\epsilon+\eta\epsilon^{2-2\ell}+\eta\epsilon^{1-2\ell})\\
&\le C_{\ell,R}(\epsilon+\eta\epsilon^{1-2\ell}).
\end{aligned}
\]

Choose \(\epsilon=\eta^{1/(2\ell)}\) and use the Kantorovich dual formula for \(W_1\) to prove \eqref{eq:source-2-1}. For \(\eta=0\), let \(\epsilon\downarrow0\). No weight appears in any constant. If \(\ell\le k\) and \(\eta\le1\), the conclusion also implies \(W_1(A,Q)\le C_{k,R}\eta^{1/(2k)}\). \end{proof}

\subsection{Deconvolution and the variance exponent}

Let \(P,Q\) have at most \(k\) atoms in \([-R,R]\), let \(s,t\in[0,b]\) be \textbf{variances}, and set

\[
\delta=\max_{1\le r\le2k}
|m_r(P*N(0,s))-m_r(Q*N(0,t))|.
\]

We prove

\begin{equation}\label{eq:source-2-6}
|s-t|\le C\delta^{1/k},\qquad
W_1(P,Q)\le C\delta^{1/(2k)}.
\end{equation}

Assume \(s\ge t\), put \(h=s-t\), and define \(A=P*N(0,h)\). The finite deconvolution polynomials are

\[
H_{r,t}(x)=r!\sum_{j=0}^{\lfloor r/2\rfloor}
\frac{(-t/2)^j x^{r-2j}}{j!(r-2j)!}.
\]

The generating identity \(\sum_rH_{r,t}(x)z^r/r!=e^{zx-tz^2/2}\), or direct finite polynomial expansion, gives
\(E H_{r,t}(Y+\sqrt t G)=E Y^r\) for an independent standard normal \(G\). Hence

\begin{equation}\label{eq:source-2-7}
\eta:=\max_{1\le r\le2k}|m_r(A)-m_r(Q)|
\le C_{k,b}\delta.
\end{equation}

The order-zero moment difference is exactly zero.

List the support points of \(Q\), repeating any of them as necessary to get a list \(q_1,\ldots,q_k\), and let \(p(x)=\prod_{j=1}^k(x-q_j)\). Its square vanishes on \(Q\) and has coefficient sum at most \((1+R)^{2k}\). Thus

\begin{equation}\label{eq:source-2-8}
E_Ap^2\le(1+R)^{2k}\eta.
\end{equation}

For any fixed \(x\), the polynomial \(z\mapsto p(x+\sqrt h z)\) has leading coefficient \(h^{k/2}\). In the monic probabilists' Hermite basis its degree-\(k\) coefficient is also \(h^{k/2}\). Hermite orthogonality therefore gives

\[
E[p(x+\sqrt h G)^2]\ge k!h^k.
\]

Average over \(x\sim P\) and combine with \eqref{eq:source-2-7}–\eqref{eq:source-2-8} to obtain the variance assertion of \eqref{eq:source-2-6}. If \(\delta\) is small enough that \(\eta\le1\), apply \eqref{eq:source-2-1} to \(A,Q\), using the actual number of distinct support points of \(Q\). Then

\[
W_1(P,Q)\le W_1(P,A)+W_1(A,Q)
\le \sqrt h E|G|+C\eta^{1/(2k)}
\le C\delta^{1/(2k)}.
\]

For the remaining discrepancy range use \(W_1(P,Q)\le2R\) and \(|s-t|\le b\), increasing the constant. If \(\delta=0\), \eqref{eq:source-2-8} first gives \(h=0\), and the zero-error case of \eqref{eq:source-2-1} gives \(P=Q\). Interchange the two models when \(t>s\). This proves \eqref{eq:source-2-6} in every case.

\subsection{Relation to the Scalar Moment Comparison}

\citet{wu2020}, Proposition 3 (author-manuscript pp.~17--18), gives precisely the exponents in \eqref{eq:source-2-6} for bounded mixing supports and bounded common standard deviations. Normalize by

\[
L=\max\{1,R,\sqrt b\},\qquad \widetilde X=X/L.
\]

Supports then lie in \([-1,1]\), standard deviations are at most one, and rescaled moment discrepancies are at most \(\delta\). Scaling back multiplies \(W_1\) by \(L\) and variance differences by \(L^2\). The proposition has no minimum weight or separation hypothesis. The displayed proposition numbering refers to arXiv:1807.07237v2. If “\(k\)-atomic” is read as exactly \(k\), add distinct atoms of total mass \(\epsilon\), scale the old weights by \(1-\epsilon\), and let \(\epsilon\downarrow0\). Compact support ensures convergence of moments and \(W_1\). Coincident projected atoms are first combined. The preceding subsections give a direct proof covering these cases, so the theorem does not depend on an interpretation of that terminology.

\section{Directional Estimation and Measurable Fitting}\label{app:estimator}

Put \(q=2k\). For every unit \(u\) and \(1\le r\le q\),

\begin{equation}\label{eq:source-4-1}
\E|u^TX|^{2r}
\le2^{2r-1}\{R^{2r}+b^r(2r-1)!!\}.
\end{equation}

This follows from \(|u^TZ|\le R\), \(u^T\varepsilon\sim N(0,u^T\Sigma u)\), and \(|x+y|^{2r}\le2^{2r-1}(|x|^{2r}+|y|^{2r})\). Let \(V=\max\{1,\max_{1\le r\le q}2^{2r-1}(R^{2r}+b^r(2r-1)!!)\}\). It bounds all relevant variances independently of \(d\).

\subsection{One fixed sphere net controls every genuine moment polynomial}

Let \(c_r=r^r/r!\), \(L_q=\max_{r\le q}r c_r\), and \(\zeta=(2L_q)^{-1}\). Fix a deterministic Euclidean \(\zeta\)-net \(\mathcal N_d\) of the unit sphere with

\begin{equation}\label{eq:source-4-2}
|\mathcal N_d|\le(1+2/\zeta)^d.
\end{equation}

A maximal \(\zeta\)-separated subset has this property: the disjoint balls of radius \(\zeta/2\) around its points fit inside the ball of radius \(1+\zeta/2\). Maximality implies the covering property.

For every real homogeneous polynomial \(h\) of degree \(r\le q\),

\begin{equation}\label{eq:source-4-3}
\sup_{\|u\|=1}|h(u)|\le2\max_{v\in\mathcal N_d}|h(v)|.
\end{equation}

Indeed, write \(h(u)=T(u,\ldots,u)\) with \(T\) symmetric \(r\)-linear, and put \(M=\sup_{\|u\|=1}|h(u)|\). The polarization identity

\[
T(x_1,\ldots,x_r)=\frac1{2^r r!}
\sum_{\epsilon\in\{-1,1\}^r}
\left(\prod_i\epsilon_i\right)h\left(\sum_i\epsilon_i x_i\right)
\]

bounds \(|T(x_1,\ldots,x_r)|\) by \(c_rM\) when all \(\|x_i\|\le1\). Telescoping then gives
\(|h(u)-h(v)|\le r c_rM\|u-v\|\) for unit \(u,v\). At a maximizer choose a net point within \(\zeta\), obtaining \(M\le\max_{\mathcal N_d}|h|+M/2\), as claimed. This argument has no ambient-dimension factor.

\subsection{One common median-of-means partition and its integrated tail}

Let \(B\) be the smallest odd integer such that

\begin{equation}\label{eq:source-4-4}
B\ge3,\qquad
B\ge\log_2q+d\log_2(1+2/\zeta).
\end{equation}

For constants depending only on \(k\), \(c_kd\le B\le C_kd\). Suppose first that \(n\ge2B\), set \(m=\lfloor n/B\rfloor\ge n/(2B)\), and divide the first \(Bm\) samples into \(B\) fixed disjoint blocks of size \(m\). Use the same blocks for every \((r,v)\in\{1,\ldots,q\}\times\mathcal N_d\). Let \(M_{r,v}\) be the median of the block averages of \((v^TX)^r\), and define

\[
\mathcal E=\max_{r,v}|M_{r,v}-\E_{P,\Sigma}(v^TX)^r|,
\qquad t_0=8\sqrt{V/m}\le8\sqrt{2VB/n}.
\]

For every \(s\ge1\),

\begin{equation}\label{eq:source-4-5}
\Pr(\mathcal E>t_0s)\le2^{-B}s^{-B}.
\end{equation}

For a fixed coordinate, Chebyshev bounds a block's error probability at \(t_0s\) by \(1/(64s^2)\). A median error requires at least \((B+1)/2\) bad independent blocks. Counting subsets by \(2^B\) bounds its probability by

\[
2^B(64s^2)^{-(B+1)/2}\le4^{-B}s^{-B}.
\]

There are at most \(q|\mathcal N_d|\le2^B\) coordinates; a union bound proves \eqref{eq:source-4-5}. Independence between coordinates is not used.

For every \(0<p\le1\), integrate this tail:

\begin{equation}\label{eq:source-4-6}
\begin{aligned}
\E\mathcal E^p
&=p\int_0^\infty t^{p-1}\Pr(\mathcal E>t)\,dt\\
&\le t_0^p\left(1+\frac{p2^{-B}}{B-p}\right)
\le2t_0^p
\le C_{k,R,b,p}(d/n)^{p/2}.
\end{aligned}
\end{equation}

The number of blocks is fixed once \(n,d,k\) are specified. It does not change with \(s\) or \(p\). Integrating this tail gives the expected-loss bound for that same fixed estimator, including its low-probability contribution. Since the logarithm of the net size is \(O_k(d)\), no \(\log n\) or \(\log d\) loss occurs.

\subsection{A single globally feasible measurable estimator}

Represent a mixing law by \(k\) labeled atoms, allowing repeated locations and zero weights. Use the compact finite-dimensional parameter space

\begin{equation}\label{eq:source-5-1}
\Theta_d=\Delta_k\times(\overline B_R)^k
\times\{S=S^T:aI_d\preceq S\preceq bI_d\},
\end{equation}

where \(\Delta_k=\{w\in[0,1]^k:\sum_jw_j=1\}\). This covers exactly the model in Section~\ref{sec:model}. For \(\theta=(w,\mu,S)\), its projected moment is the continuous polynomial

\begin{equation}\label{eq:source-5-2}
m_{r,v}(\theta)=\sum_{j=1}^k w_j
\sum_{h=0}^{\lfloor r/2\rfloor}
{r\choose2h}(2h-1)!!(v^TSv)^h(v^T\mu_j)^{r-2h},
\end{equation}

with \((-1)!!=1\).

Fix a countable dense sequence \((\theta_j)_{j\ge1}\) in \(\Theta_d\). Concretely use rational simplex weights, rational locations strictly inside the ball, and rational symmetric matrices with \(aI_d\prec S\prec bI_d\). Density at the boundary follows by first contracting locations toward zero and moving covariance matrices toward an interior scalar matrix, and then approximating by rationals inside the constraints. This works even when \(R,a,b\) are irrational.

For \(n\ge2B\) and the targets of Appendix~\ref{app:estimator} define

\[
A_j=\max_{r,v}|m_{r,v}(\theta_j)-M_{r,v}|,
\quad A_*=\inf_{j\ge1}A_j,
\quad\tau_n=\sqrt{B/n}>0,
\]

\begin{equation}\label{eq:source-5-3}
J=\min\{j:A_j\le A_*+\tau_n\},
\qquad (\widehat P,\widehat\Sigma)=(P_{\theta_J},S_{\theta_J}).
\end{equation}

For fixed data the fitting objective is continuous on the compact set \(\Theta_d\). The dense infimum equals the full minimum and is at most \(\mathcal E\), since the true parameter is feasible. Positive tolerance ensures that \(J\) exists. Each block mean, median, and \(A_j\) is Borel measurable in the data. So is the countable infimum \(A_*\), and

\[
\{J=j\}=\{A_j\le A_*+\tau_n\}
\cap\bigcap_{i<j}\{A_i>A_*+\tau_n\}
\]

is measurable. The countably valued joint estimator is therefore measurable and belongs to the original parameter class for every sample. No independently fitted projections are substituted for a joint model.

On the net,

\[
|m_{r,v}(\theta_J)-m_{r,v}(\theta_{\rm true})|
\le A_J+\mathcal E\le2\mathcal E+\tau_n.
\]

For each \(r\), the difference between these two model moments is a homogeneous polynomial of degree \(r\) in \(v\). Apply \eqref{eq:source-4-3}, after fitting, to get

\begin{equation}\label{eq:source-5-4}
D_{2k}(F_{\widehat P,\widehat\Sigma},F_{P,\Sigma})
\le4\mathcal E+2\tau_n.
\end{equation}

Polynomial norming has been applied to the fitted model moment difference. For \(0<p\le1\), subadditivity of \(x^p\), \eqref{eq:source-4-6}, and \(B\le C_kd\) imply

\begin{equation}\label{eq:source-5-5}
\E D_{2k}^p\le C_{k,R,b,p}(d/n)^{p/2}.
\end{equation}

Apply \eqref{eq:source-1-3} and \eqref{eq:source-5-5} with \(p=1/(2k)\) and \(p=1/k\) to obtain

\begin{equation}\label{eq:source-5-6}
\E W_1(\widehat P,P)\le C(d/n)^{1/(4k)},\qquad
\E\|\widehat\Sigma-\Sigma\|_{\rm op}\le C(d/n)^{1/(2k)}
\quad(n\ge2B).
\end{equation}

For \(n<2B\), define the estimator to be \((\delta_0,aI_d)\). Its losses are at most \(R\) and \(b-a\). In this branch \(n<2C_kd\), so for either exponent \(\alpha\) in \eqref{eq:source-1-1}–\eqref{eq:source-1-2},

\[
\min\{1,(d/n)^\alpha\}\ge\min\{1,(2C_k)^{-\alpha}\}>0.
\]

Thus bounded losses give the claimed saturated upper bounds also in this branch. In the other branch \(n\ge2B\ge2d\), so \eqref{eq:source-5-6} already has the unsaturated target. This proves both upper bounds for every \(n,d\ge1\).

\section{Positive Quadrature and Gaussian Information}\label{app:quadrature}

\subsection{Positive Gaussian quadrature for every fixed \texorpdfstring{$k$}{k}}

Let \(H_j\) be the monic probabilists' Hermite polynomials, so

\[
e^{zx-z^2/2}=\sum_{j\ge0}H_j(x)z^j/j!,\qquad
E H_i(G)H_j(G)=j!\,1_{i=j},\quad G\sim N(0,1).
\]

The degree-$k$ polynomial \(H_k\) has $k$ distinct real roots \(a_1,\ldots,a_k\). To see this, if it had fewer than $k$ distinct sign-changing real roots, multiply their linear factors to form a polynomial \(q\) of degree less than $k$. The real polynomial \(H_kq\) has a constant nonzero sign except at finitely many zeros, contradicting \(E[H_k(G)q(G)]=0\). Thus it has $k$ sign-changing roots and all are simple.

Let \(l_i\) be the Lagrange polynomial of degree \(k-1\) with \(l_i(a_j)=1_{i=j}\), and set

\begin{equation}\label{eq:source-6-4}
w_i=E l_i(G)=E l_i(G)^2>0.
\end{equation}

The equality holds because \(l_i^2-l_i\) is divisible by \(H_k\) with quotient of degree at most \(k-2\). For any polynomial \(f\) of degree at most \(2k-1\), division by \(H_k\), orthogonality, and interpolation show that

\begin{equation}\label{eq:source-6-5}
E f(G)=\sum_{i=1}^k w_i f(a_i).
\end{equation}

In particular \(\sum_iw_i=1\). Let \(A\) have this positive quadrature law, and let \(B_k=\max_i|a_i|\). Then \(B_k\ge1\), \(EA=0\), \(EA^2=1\), and

\begin{equation}\label{eq:source-6-6}
E A^j=E G^j\quad(0\le j\le2k-1),\qquad
E|A|\ge B_k^{-1}.
\end{equation}

The latter follows from \(A^2\le B_k|A|\). Parity of \(H_k\) and uniqueness of the interpolation weights imply symmetry of \(A\). For \(k=2\), \(H_2(x)=x^2-1\) and this is the law in \eqref{eq:source-6-1}.

The first unmatched Hermite moment is nonzero and explicit:

\begin{equation}\label{eq:source-6-7}
E H_{2k}(A)=-k!.
\end{equation}

Indeed, \(H_k^2-H_{2k}\) has degree at most \(2k-1\) and Gaussian expectation \(k!\), whereas \(H_k(A)=0\). This also says \(EA^{2k}-EG^{2k}=-k!\).

Fix the positive constant

\[
t_{\max}=\min\{1,R/B_k,\sqrt{(v-a)/2},\sqrt v/2\}.
\]

For \(0<t\le t_{\max}\) and any unit \(u\), define

\begin{equation}\label{eq:source-6-8}
P_{u,t}=\mathcal L(tAu),\qquad
\Sigma_{u,t}=vI_d-t^2uu^T.
\end{equation}

The law has exactly $k$ distinct bounded atoms with positive weights. Its noise covariance has the same admissible eigenvalues as in Section~\ref{sec:lower}. The observation density is a positive normalized finite Gaussian mixture. Its projection along \(u\) is
\(tA+\sqrt{v-t^2}G\); the orthogonal part is an independent \(N(0,vI_{u^\perp})\).

Replacing \(A\) by an independent standard normal in a polynomial of total degree at most \(2k-1\) does not change its expectation, by \eqref{eq:source-6-6}. After that replacement the whole vector has law \(Q_0\). Hence all raw multivariate moments through total degree \(2k-1\) agree exactly with \(Q_0\), and for every \(z\),

\begin{equation}\label{eq:source-6-9}
E(z^TX)^{2k}-E_{Q_0}(z^TX)^{2k}
=-k!t^{2k}(z^Tu)^{2k}.
\end{equation}

In particular the compensated total covariance is always \(vI_d\).

\subsection{Likelihood integrability, tails, and an exact nonnegative remainder}

Put \(s=t/\sqrt v\le1/2\), and let \(q_s\) be the density of \(Y_s=sA+\sqrt{1-s^2}G\). Relative to the standard normal density \(\phi\), its likelihood ratio is

\begin{equation}\label{eq:source-6-10}
L_s(x)=E K_s(x,A),\qquad
K_s(x,z)=\frac1{\sqrt{1-s^2}}
\exp\left\{\frac{-s^2x^2+2sxz-s^2z^2}{2(1-s^2)}\right\}.
\end{equation}

Completing the square, for \(s>0\), gives

\[
\frac{-s^2x^2+2sxz-s^2z^2}{2(1-s^2)}
=-\frac{s^2}{2(1-s^2)}(x-z/s)^2+\frac{z^2}{2}.
\]

Consequently, including \(s=0\) by its direct value, for every \(s\in[0,1/2]\) and every \(x\),

\begin{equation}\label{eq:source-6-11}
0<L_s(x)\le C_0:=\frac2{\sqrt3}e^{B_k^2/2}.
\end{equation}

In particular \(L_s\in L^2(\phi)\), uniformly, \(\int L_s\phi=1\), and the explicit Gaussian tail bound

\begin{equation}\label{eq:source-6-12}
\int_{|x|>T}L_s(x)^2\phi(x)\,dx
\le C_0^2\Pr(|G|>T)\le2C_0^2e^{-T^2/2}
\quad(T\ge0)
\end{equation}

controls the entire likelihood tail.

For every integer \(m\ge0\), a finite polynomial expansion of the Gaussian generating identity yields

\begin{equation}\label{eq:source-6-13}
\int H_m(x)K_s(x,z)\phi(x)\,dx=s^mH_m(z).
\end{equation}

For clarity, integrate \(e^{wx-w^2/2}\) under \(N(sz,1-s^2)\), obtaining \(e^{swz-s^2w^2/2}\), and take the \(m\)-th derivative at zero. Gaussian exponential moments justify this finite differentiation.

Completeness and orthogonality of the Hermite polynomials in Gaussian \(L^2\), with \eqref{eq:source-6-11}, now give the convergent expansion

\[
L_s=\sum_{m\ge0}\frac{s^m E H_m(A)}{m!}H_m
\quad\text{in }L^2(\phi).
\]

One direct completeness justification is that if \(f\in L^2(\phi)\) is orthogonal to every polynomial, then \(\int f(x)\phi(x)e^{zx}\,dx\) is an entire function: Cauchy–Schwarz and Gaussian exponential integrability give convergence locally uniformly in complex \(z\), also for all derivatives. Every derivative at zero vanishes. The entire function is zero; restricting to the imaginary axis and using uniqueness of Fourier transforms shows \(f\phi=0\).

Since \(EH_m(A)=0\) for \(1\le m\le2k-1\), Parseval supplies the \textbf{exact} nonnegative remainder

\begin{equation}\label{eq:source-6-14}
\chi^2(q_s,\phi)=\sum_{m\ge2k}
\frac{s^{2m}(EH_m(A))^2}{m!}.
\end{equation}

At \(s_0=1/2\), \(\|L_{s_0}\|_2^2\le C_0\int L_{s_0}\phi=C_0\). Therefore, for \(0\le s\le s_0\), each nonnegative summand in \eqref{eq:source-6-14} is at most \((s/s_0)^{4k}\) times its value at \(s_0\), and

\begin{equation}\label{eq:source-6-15}
\frac{(k!)^2}{(2k)!}s^{4k}
\le\chi^2(q_s,\phi)
\le2^{4k}C_0s^{4k}.
\end{equation}

The lower inequality uses just the \(m=2k\) term and \eqref{eq:source-6-7}. The likelihood bound and the nonnegative series control the full remainder.

The observation law and the reference share the independent orthogonal Gaussian factor. A rotation and the change of scale \(\sqrt v\) leave the divergence unchanged. It follows that, uniformly in \(d,u,t\) as above,

\begin{equation}\label{eq:source-6-16}
\operatorname{KL}(F_{u,t}\|Q_0)
\le\chi^2(F_{u,t},Q_0)
\le C_\chi t^{4k},\qquad
C_\chi=2^{4k}C_0v^{-2k}.
\end{equation}

The inequality \(\operatorname{KL}\le\chi^2\) follows from \(\log x\le x-1\) and normalization. For \(k=2\), \eqref{eq:source-6-15} begins with \(s^8/6\) and proves the required \(Ct^8\) bound for \eqref{eq:source-6-2}.

Section~\ref{sec:sharpness} uses the exact moment identity \eqref{eq:source-6-9} to prove Corollary~\ref{cor:sharp}.

\section{Projective Packing and Testing}\label{app:testing}

\subsection{Distances between two compensated directions}

For unit \(u,u'\), let \(\rho=u^Tu'\). Every target atom for \(P_{u',t}\) lies on \(\operatorname{span}(u')\). Any coupling therefore incurs at least the distance of each source point to that line, giving

\begin{equation}\label{eq:source-7-1}
W_1(P_{u,t},P_{u',t})
\ge tE|A|\sqrt{1-\rho^2}
\ge(t/B_k)\sqrt{1-\rho^2}.
\end{equation}

This remains true when the quadrature law has a zero node. A direct calculation in \(\operatorname{span}(u,u')\) gives the nonzero eigenvalues of \(uu^T-u'u'^T\) as \(\pm\sqrt{1-\rho^2}\); in the collinear case the matrix is zero. Thus

\begin{equation}\label{eq:source-7-2}
\|\Sigma_{u,t}-\Sigma_{u',t}\|_{\rm op}
=t^2\sqrt{1-\rho^2}.
\end{equation}

\subsection{An exponential projective packing}

Set \(h(p)=-p\log p-(1-p)\log(1-p)\),
\(\gamma=\log2-h(1/4)>0\). Greedily select from the normalized sign vectors \(\{\pm1/\sqrt d\}^d\), removing vectors with absolute inner product greater than \(1/2\) with each selected vector. At each step at most

\[
2\sum_{j\le\lfloor d/4\rfloor}{d\choose j}
\le2e^{d h(1/4)}
\]

vectors are removed. For the last bound, if \(J\sim\operatorname{Bin}(d,1/2)\), the exponential Markov inequality at \(\lambda=\log3\) gives
\(\Pr(J\le d/4)\le e^{d h(1/4)-d\log2}\). The inner product of two sign vectors is \(1-2H/d\), where \(H\) is their Hamming distance, which explains both removed caps.

We obtain \(M\ge\tfrac12e^{\gamma d}\) directions with
\(|u_i^Tu_j|\le1/2\) for \(i\ne j\). This removes both signs and hence packs directions modulo sign, consistent with the symmetry \(P_{-u,t}=P_{u,t}\). Choose a universal integer \(d_0\) sufficiently large that, for all \(d\ge d_0\),

\begin{equation}\label{eq:source-7-3}
\log M\ge(\gamma/2)d=:c_0d,\qquad
\log M\ge4\log2.
\end{equation}

\subsection{Fano, with bounded-loss saturation included}

Choose a fixed \(\eta_0>0\) with
\(\eta_0\le t_{\max}\) and \(C_\chi\eta_0^{4k}\le c_0/8\). For \(d\ge d_0\) use the packed parameters at

\begin{equation}\label{eq:source-7-4}
t=\eta_0\min\{1,(d/n)^{1/(4k)}\}.
\end{equation}

Product KL additivity and \eqref{eq:source-6-16} give

\begin{equation}\label{eq:source-7-5}
\frac1M\sum_i\operatorname{KL}(F_{u_i,t}^{\otimes n}\|Q_0^{\otimes n})
\le nC_\chi t^{4k}
\le C_\chi\eta_0^{4k}d\le\tfrac18\log M.
\end{equation}

For completeness, if the index is uniform, its mutual information with the data is at most the average KL to any fixed reference law. This follows from the exact identity
\(M^{-1}\sum_i\operatorname{KL}(F_i\|Q)=I+\operatorname{KL}(\overline F\|Q)\).
Fano's entropy inequality then gives error probability at least
\(1-(I+\log2)/\log M\ge5/8\) for every decoder. This form uses the reference law only for bounding information; it need not be a member of the packing.

For either parameter loss, decode an arbitrary estimator to its closest packing member, resolving finite ties by index. If its loss to the true member is less than half the minimum pairwise separation, it decodes correctly by the triangle inequality. Equations \eqref{eq:source-7-1}–\eqref{eq:source-7-2} and the testing bound imply

\begin{equation}\label{eq:source-7-6}
\sup_{\mathcal M_{k,d}}E W_1(\widehat P,P)\ge c t,\qquad
\sup_{\mathcal M_{k,d}}E\|\widehat\Sigma-\Sigma\|_{\rm op}\ge c t^2.
\end{equation}

This argument applies separately to arbitrary estimators in the two minimax definitions, without requiring their outputs to lie in the model. If a transport loss is infinite the lower bound already holds. When \(n\le d\), the chosen \(t=\eta_0\) is constant and \eqref{eq:source-7-5} still holds, proving constant-order lower bounds in the saturated regime.

\subsection{Bounded Dimensions, Including One Dimension}

For \(1\le d<d_0\), compare the two admissible parameters

\[
\theta_0=(\delta_0,vI_d),\qquad
\theta_1=(P_{e_1,t},vI_d-t^2e_1e_1^T),
\quad t=\eta_1n^{-1/(4k)},
\]

where \(0<\eta_1\le t_{\max}\) and \(C_\chi\eta_1^{4k}\le1/8\). Their product KL is at most \(1/8\), so Pinsker's inequality bounds their product total variation by \(1/4\). Their parameter distances are

\[
W_1(\delta_0,P_{e_1,t})=tE|A|\ge t/B_k,
\qquad\|vI_d-\Sigma_{e_1,t}\|_{\rm op}=t^2.
\]

For any estimator, nearest-member testing and the triangle inequality give maximum risk at least one quarter of the parameter distance times \(1-\mathrm{TV}\). Hence the risks are at least constant multiples of \(n^{-1/(4k)}\) and \(n^{-1/(2k)}\). For either positive exponent \(\alpha\),

\[
\min\{1,(d/n)^\alpha\}\le d_0^\alpha n^{-\alpha}
\quad(1\le d<d_0).
\]

Thus these two-point bounds give the same saturated targets up to a fixed constant in all remaining dimensions and sample-size regimes. Combining them with \eqref{eq:source-7-6} proves both lower inequalities in \eqref{eq:source-1-1}–\eqref{eq:source-1-2}. Together with Appendix~\ref{app:estimator}, this proves Theorem~\ref{thm:minimax}. 

\bibliography{references}
\end{document}